\documentclass[10pt]{article}

\usepackage{graphicx}                           
\usepackage{color}                              
\usepackage{indentfirst}                        
\usepackage{amsmath,amssymb,amsfonts,amsthm,bm} 
\usepackage{cases}                              
\usepackage{cite}
\usepackage{mathrsfs}

\newtheorem{theorem}{Theorem}[section]
\newtheorem{definition}{Definition}[section]
\newtheorem{lemma}{Lemma}[section]

\newtheorem{example}{Example}[section]
\newtheorem{remark}{Remark}[section]

\theoremstyle{remark}
\theoremstyle{remark}
\begin{document}
\title{Simulative Suzuki-Gerghaty type contraction with $\mathcal{C}$-class functions and applications}
\author{Abdullah Eqal Al-Mazrooei \thanks{Department of Mathematics, University of Jeddah, P.O.Box 80327, Jeddah 21589, Saudi Arabia, Email: aealmazrooei@uj.edu.sa,} Azhar Hussain\thanks{
Department of Mathematics, University of Sargodha, Sargodha-40100, Pakistan. Email: hafiziqbal30@yahoo.com},\\ Muhammad Ishfaq\thanks{
Department of Mathematics, University of Sargodha, Sargodha-40100, Pakistan. Email: ishfaqahmad5632@gmail.com,} Jamshaid Ahmad\thanks{Department of Mathematics, University of Jeddah,, P.O.Box 80327, Jeddah 21589,
Saudi Arabia. Email: jkhan@uj.edu.sa.}}
\date{}
\maketitle

\noindent\textbf{ Abstract}: The aim of this paper is to introduce the notion of a Suzuki-Gerghaty type contractive mapping via simulation function along with $\mathcal{C}$-class functions and prove the existence of fixed point result. An example is given to show the validity of our results given herein. Moreover, we prove the existence of solution of nonlinear Hammerstein integral equation.  \\

\noindent{\textbf{Mathematics Subject Classification}}: 54H25, 47H10\\

\noindent{\textbf{ Keywords}}: Simulation functions, $\mathcal{C}$-class function, nonlinear integral equation.

\section{Introduction and Literature Review}
In 2015, Khojasteh {\it et al.} \cite{stojan2} introduced a function $\zeta : [0, \infty) \times  [0, \infty) \rightarrow \mathbb{R}$, satisfying the following assertions:
\begin{description}
  \item[$(\zeta_{1})$] $\zeta(0, 0)=0$;
  \item[$(\zeta_{2})$]  $\zeta(t, s)<s-t$ for all $t, s>0$;
  \item[$(\zeta_{3})$]  If $\{t_{n}\}, \{s_{n}\}$ are sequences in $(0, \infty)$ such that $\lim\limits_{n\rightarrow \infty}t_{n}=\lim\limits_{n\to \infty}s_{n}>0$ then $$\lim\limits_{n\rightarrow \infty}sup~\zeta(t_{n}, s_{n})<0$$
\end{description}
called the simulation function and define the notion of $\mathcal{Z}$-contraction with respect to the function $\zeta$ to generalized Banach contraction principle \cite{banach} and unify several known contractions involving the combination of $d(Tx, Ty)$ and $d(x, y)$.
They presented the following contractive condition via simulation function and prove a fixed point theorem:
\begin{definition}\cite{stojan2}
Let $(X, d)$ be a metric space, $T: X \rightarrow X$ a mapping and $\zeta$ a simulation function. Then $T$ is called a $\mathcal{Z}$-contraction with respect to $\zeta$ if it satisfies
\begin{equation}\label{Eq01}
\zeta(d(Tu, Tv), d(u, v)) \geq 0 ~~for~ all~ u, v \in X.
\end{equation}
\end{definition}
\begin{theorem}\label{T1}\cite{stojan2}
Let $(X, d)$ be a complete metric space and $T : X \rightarrow X$ be a $\mathcal{Z}$-contraction with respect to $\zeta$. Then $T$ has a unique fixed point $u\in X$ and for every $x_{0} \in X$, the Picard sequence $\{x_{n}\}$ where $x_{n}=Tx_{n-1}$ for all $n \in \mathbb{N}$ converges to this fixed point of $T$.
\end{theorem}
\begin{example}\cite{stojan2}
Let $ \zeta_{i} : [0, \infty) \times [0, \infty)\rightarrow \mathbb{R},~~i=1, 2, 3$ be defined by
\begin{description}
  \item[$(i)$] $\zeta_{1}(t, s)=\lambda s-t$, where $\lambda \in (0, 1)$;
  \item[$(ii)$]  $\zeta_{2}(t, s)= s\varphi(s)-t$, where $\varphi:[0, \infty)\rightarrow [0, 1)$ is a mapping such that $\lim\limits_{t\rightarrow r^{+}}sup~\psi(t)<1$ for all $r>0$;
  \item[$(iii)$]  $\zeta_{3}=s-\psi(s)-t$, where $\psi:[0, \infty)\rightarrow [0, \infty)$ is a continuous function such that $\psi(t)=0$ if and only if $t=0.$
\end{description}
Then $\zeta_{i}$ for $i=1, 2, 3$ are simulation functions.
\end{example}
Karapinar \cite{Erdal} presented some fixed point results in the setting of a complete metric spaces by defining a new contractive condition via admissible mapping imbedded in simulation function. Hakan \cite{HM} {\it et al.} introduced the generalized simulation function on a quasi metric space and presented a fixed point theorem.
Rold$\acute{a}$n-L$\acute{o}$pez-de-Hierro {\it et al.} \cite{Roldan} modified the notion of a simulation function by replacing $(\zeta_{3})$ by $(\zeta'_{3})$,
\begin{description}
\item[($\zeta'_{3}$)]: if $\{t_{n}\}, \{s_{n}\}$ are sequences in $(0, \infty)$ such that $\lim_{n\rightarrow \infty}t_{n}=\lim\limits_{n\rightarrow \infty}s_{n}>0$ and $t_{n}<s_{n}$, then $$\lim\limits_{n\rightarrow \infty}sup\zeta(t_{n}, s_{n})<0.$$
\end{description}
The function $\zeta: [0, \infty) \times [0, \infty) \rightarrow \mathbb{R}$ satisfying $(\zeta_{1}-\zeta{2})$ and $(\zeta'_{3})$ is called simulation function in the sense of Rold$\acute{a}$n-L$\acute{o}$pez-de-Hierro.

\begin{definition}\cite{arslan}\label{D4.2}
A mapping $\mathcal{G}:[0,+\infty)^{2}\rightarrow \mathbb{R}$ is called a
$\mathcal{C}$-class function if it is continuous and satisfies the following
conditions:
\begin{description}
  \item [(1)]  $\mathcal{G}(s, t)\leq s$;
  \item [(2)]  $\mathcal{G}(s, t)= s$ implies that either $s=0$ or $t=0$,
      for all $s, t \in [0, +\infty)$.
\end{description}
\end{definition}
\begin{definition}\cite{arslan2}\label{D4.3}
A mapping $\mathcal{G}:[0,+\infty)^{2}\rightarrow\mathbb{R}$ has the
property $\mathcal{C}_{\mathcal{G}}$, if there exists and $\mathcal{C}_{\mathcal{G}}\geq 0$ such
that
\begin{description}
  \item [(1)]  $\mathcal{G}(s, t)> \mathcal{C}_{\mathcal{G}}$ implies $s>t$;
  \item [(2)]  $\mathcal{G}(s, t)\leq \mathcal{C}_{\mathcal{G}}$, for all $t\in[0,
      +\infty)$.
\end{description}
\end{definition}
Some examples of $\mathcal{C}$-class functions that have property $\mathcal{C}_{\mathcal{G}}$ are
as follows:
\begin{description}
  \item [(a)] $\mathcal{G}(s, t)= s-t,~ \mathcal{C}_{\mathcal{G}}=r, r\in[0,
      +\infty)$;
  \item [(b)] $\mathcal{G}(s, t)= s-\frac{(2+t)t}{(1+t)},
      \mathcal{C}_{\mathcal{G}}=0$;
  \item[(c)] $\mathcal{G}(s, t)= \frac{s}{1+kt}, k\geq1,
      \mathcal{C}_{\mathcal{G}}=\frac{r}{1+k}, r\geq 2$.
\end{description}
For more examples of $\mathcal{C}$-class functions that have property $\mathcal{C}_{\mathcal{G}}$
see \cite{arslan3, ZM, arslan2}.
\begin{definition}\cite{arslan2}\label{D4.5}
A $\mathcal{C}_{\mathcal{G}}$ simulation function is a mapping
$\mathcal{G}:[0,+\infty)^{2}\rightarrow\mathcal{R}$ satisfying the following
conditions:
\begin{description}
  \item [(1)]  $\zeta(t, s)< \mathcal{G}(s, t)$ for all $t, s>0$, where
      $\mathcal{G}:[0, +\infty)^{2}\rightarrow \mathbb{R}$ is a $\mathcal{C}$-class
      function;
  \item [(2)] if $\{t_{n}\}, \{s_{n}\}$ are sequences in $(0, +\infty)$
      such that $\lim\limits_{n\rightarrow \infty}t_{n}=
      \lim\limits_{n\rightarrow \infty}s_{n}>0 $, and $t_{n}<s_{n}$, then
      $\lim\limits_{n\rightarrow\infty}\sup\zeta(t_{n},
      s_{n})<\mathcal{C}_{\mathcal{G}}$.
\end{description}
\end{definition}
Some examples of simulation functions and $\mathcal{C}_{\mathcal{G}}$-simulation
functions are:
\begin{description}
  \item [(1)]  $\zeta(t, s)= \frac{s}{s+1}-t$ for all $t, s>0$.
  \item [(2)] $\zeta(t, s)= s-\phi(s)-t$ for all $t, s>0$, where
      $\phi:[0, +\infty)\rightarrow [0, +\infty)$ is a lower semi
      continuous function and $\phi(t)=0$ if and only if $t=0$.
\end{description}
For more examples of simulation functions and $\mathcal{C}_{\mathcal{G}}$-simulation
functions see \cite{arslan3, Roldan, stojan2, arslan2, vetro, wang}.

We denote by $\mathcal{F}$ the class of all functions $ \beta : [0,\infty) \rightarrow [0,1)$ satisfying
$ \beta (t_{n}) \to 1$, implies $t_{n} \rightarrow 0$ as $n \to \infty$.

\begin{definition}\cite{ger}
Let $(X,d)$ be a metric space. A map $T:X\rightarrow X$ is called Geraghty contraction if there exists $\beta\in \mathcal{F}$ such that for all $x, y\in X$,
$$d(Tx,Ty)\leq\beta(d(x,y))d(x,y).$$
\end{definition}
By using such maps Geraghty {\it et al.} \cite{ger} proved the following fixed point result:
\begin{theorem}
Let $(X,d)$ be a complete metric space. Mapping $T:X \rightarrow X$ is Geraghty contraction.
Then $T$ has a fixed point $x\in X$, and $\{T^{n}x_{1}\}$ converges to $x$.
\end{theorem}
\begin{definition}\cite{samet}
Let $T:A\rightarrow B$ be a map and $\alpha :X \times X \rightarrow \mathbb{R}$ be a function. Then $T$ is said to be $\alpha$-admissible if $\alpha(x,y)\geq 1$
 implies $\alpha(Tx,Ty)\geq 1.$
\end{definition}
\begin{definition}\cite{karap}
An $\alpha$-admissible map $T$ is said to be triangular $\alpha$-admissible if $\alpha (x,z)\geq 1$ and $\alpha (z,y)\geq 1$  implies $\alpha(x,y)\geq 1$
\end{definition}

Cho et al. \cite{cho} generalized the concept of Geraghty contraction to $\alpha$-Geraghty contraction and prove the fixed point theorem for such contraction.
\begin{definition}\cite{cho}
Let $(X,d)$ be a metric space, and let $\alpha : X \times X \rightarrow \mathbb{R}$ be a function. A map $T:X\rightarrow X$ is called $\alpha$-Geraghty contraction if there exists $\beta\in \mathcal{F}$ such that for all $x, y\in X$,
$$\alpha (x,y)d(Tx,Ty)\leq\beta(d(x,y))d(x,y).$$
\end{definition}
\begin{theorem}\cite{cho}
Let $(X,d)$ be a complete metric space, $\alpha: X\times X \rightarrow \mathbb{R}$ be a function. Define a map $T:X \rightarrow X$ satisfying the following conditions:
\begin{enumerate}
  \item $T$ is continuous $\alpha$-Geraghty contraction;
  \item $T$ be a triangular $\alpha$-admissible;
  \item there exists $x_{1}\in X$ such that $\alpha(x_{1}, Tx_{1})\geq 1$;
\end{enumerate}
Then $T$ has a fixed point $x\in X$, and $\{T^{n}x_{1}\}$ converges to $x$.
\end{theorem}
On the other hand Suzuki \cite{SUZUKI1} generalized the Banach contraction principle by using the following contractive condition:
\begin{definition}
Let $(X, d)$ be a metric space, the mapping $T: X \rightarrow X$ is called Suzuki contraction if for all $x, y\in X$
\begin{equation}
\frac{1}{2}d(x, Tx)<d(x, y) ~~ \text{implies} ~~ d(Tx, Ty)<d(x, y).
\end{equation}
\end{definition}
In 2017, Kumam {\it et al.} \cite{kumam}, following Suzuki \cite{SUZUKI1}, introduced the notion Suzuki-type $\mathcal{Z}$-contraction and gave the following generalization:
\begin{definition} \cite{kumam}
Let $T:X\rightarrow X$ be a mapping and $x_{0} \in X$ be arbitrary. Then $T$ is said to be possess property $(K)$ if for a bounded Picard sequence $\{x_{n}=Tx_{n-1};n \in \mathbb{N}\}$ there exists subsequence $\{x_{mk}\}$ and $\{x_{nk}\}$ of $\{x_{n}\}$ such that $\lim\limits_{n\rightarrow \infty}d(x_{mk}, x_{nk})=C>0$ where $m_{k}>n_{k}>k$, $k \in \mathbb{N}$ then
$$\frac{1}{2}d(x_{mk-1}, x_{mk})<d(x_{mk-1}, x_{nk-1})$$
holds.
\end{definition}
\begin{definition} \cite{kumam}
Let $(X, d)$ be a metric space, $T: X \rightarrow X$ a mapping and $\zeta \in \mathcal{Z}$. Then $T$ is called a Suzuki type $\mathcal{Z}$-contraction with respect to $\zeta$ if the following condition is satisfied
\begin{equation}\label{p}
\frac{1}{2}d(x, Tx) < d(x, y) \Rightarrow \zeta(d(Tx, Ty), d(x, y)) \geq 0
\end{equation}
for all $x, y \in X$, with $x \neq y$
\end{definition}
\begin{theorem}\cite{kumam}
Let $(X, d)$ be a complete metric space and $T:X\rightarrow X$ be a Suzuki type $\mathcal{Z}$-contraction with respect to $\zeta$. Then $T$ has a unique fixed point $x \in X$ and for every $x_{0} \in X$ the Picard sequence $\{x_{n}\}$, where $x_{n}=Tx_{n-1}$ for all $n \in \mathbb{N}$ converges to the fixed point of $T$, provided that $F$ has property $(K)$.
\end{theorem}
\begin{lemma}\cite{stojan}\label{L4.2}
Let $(X, d)$ be a metric space and let $\{x_{n}\}$ be a sequence in $X$ such
that
\begin{equation}
\lim\limits_{n\rightarrow\infty}d(x_{n}, x_{n+1})=0.
\end{equation}
If $\{x_{n}\}$ is not a Cauchy sequence in $X$, then there exists
$\varepsilon >0$ and two sequences $x_{m(k)}$ and $x_{n(k)}$ of positive
integers such that $x_{n(k)}>x_{m(k)}>k$ and the following sequences tend to
$\varepsilon^{+}$ when $k\rightarrow\infty$:
$$d(x_{m(k)},x_{n(k)}), d(x_{m(k)},x_{n(k)+1}),d(x_{m(k)-1},x_{n(k)}),$$ $$d(x_{m(k)-1},x_{n(k)+1}),d(x_{m(k)+1},x_{n(k)+1}).$$
\end{lemma}
The purpose of this paper is to generalized the results of Padcharoen {\it et al.} \cite{poom} in the context of simulative $\mathcal{C}$-class functions with Geraghty's view.
\section{Main Results}
We begin with the following notion:
\begin{definition}\label{D4.1}
Let $(X, d)$ be a metric space and $T:X\rightarrow X$ be a mapping.
$\alpha: X\times X\rightarrow [0, \infty)$  be a function and $\zeta\in
\mathcal{Z}$. Mapping $T$ is called Suzuki type $\mathcal{Z}_{(\alpha,
\mathcal{G})}$-Geraghty contraction with respect to $\zeta$ if for all distinct $x, y\in X$
\begin{equation}\label{E4.1}
\frac{1}{2}d(x, Tx)<d(x, y) ~\text{implies} ~ \zeta(\alpha(x, y)d(Tx, Ty), \beta(M(x, y)) M(x, y))\geq \mathcal{C}_{\mathcal{G}},
\end{equation}
where
$$ M(x, y)= \max\{d(x, y), d(x, Tx), d(y, Ty)\}.$$
\end{definition}
\begin{remark}
By definition of $\mathcal{C}_{\mathcal{G}}$-simulation function with $\mathcal{C}_{\mathcal{G}}=0$ and $\mathcal{G}(s, t)=s-t$, we can have from \eqref{E4.1} that
\begin{equation}\label{E4.1'}
\frac{1}{2}d(x, Tx)<d(x, y) ~\text{implies} ~ \alpha(x, y)d(Tx, Ty)< \beta(M(x, y)) M(x, y)).
\end{equation}
\end{remark}
Following lemma will be needed in the sequel:
\begin{lemma}\label{L3.1}
Let $T:X\rightarrow X$ be a triangular $\alpha$-admissible map. Assume that
there exists $x_{1} \subset X$ such that $ \alpha(x_{1}, Tx_{1})\geq 1
$.Define a sequence $\{x_{n}\}$ by $x_{n+1}= Tx_{n}$. Then we have
$\alpha(x_{n}, x_{m})\geq 1$ for all $m, n \in \mathbb{N}$ with $n < m$.
\begin{proof}
Since there exist $x_{1}\in X$ such that $\alpha(x_{1}, Tx_{1})\geq 1$, then by the
$\alpha$-admissibility of $T$, we have $\alpha(Tx_{1}, T^{2}x_{1})=\alpha(x_{2}, x_{3})\geq 1$. By continuing this process, we get
$\alpha(x_{n}, x_{n+1})\geq 1$ for all $n\in\mathbb N\cup\{0\}$. Suppose that $n < m$. Since $\alpha(x_{n}, x_{n+1})\geq 1$
and $\alpha(x_{n+1}, x_{n+2})\geq 1$, then using the triangular $\alpha$-admissibility of $T$ we have
$\alpha(x_{n}, x_{n+2}) \geq 1$. Again, since $\alpha(x_{n}, x_{n+2})\geq 1$ and $\alpha(x_{n+2}, x_{n+3})\geq 1$, then we
deduce $\alpha(x_{n}, x_{n+3})\geq 1$. By continuing this process, we get $\alpha(x_{n}, x_{m})\geq 1$.
\end{proof}
\end{lemma}

To prove the existence of fixed point result, we need the following lemmas:
\begin{lemma}\label{L4.1}
Let $(X, d)$ be a complete metric space, $\alpha:X\times X\rightarrow [0,
\infty)$ and $T: X \rightarrow X$ be two  functions. Suppose that the
following conditions are satisfied:
\begin{description}
  \item [(1)]  $T$ is Suzuki type $\mathcal{Z}_{(\alpha,
      \mathcal{G})}$-Geraghty contraction with respect to $\zeta$;
  \item [(2)] $T$ is triangular $\alpha$-admissible;
  \item[(3)] for all $x,y\in Fix(T)$, there exists $z\in X$ such that $\alpha(x, z)\geq 1$ and
$\alpha(y, z)\geq 1$.
\end{description}
Then $T$ has at most one fixed point.
\begin{proof}
Suppose that $x, y \in Fix(T)$, with $x\neq y$.  Since
$$0=\frac{1}{2}d(x, Tx)<d(x, y),$$
then by (\ref{E4.1}), we obtain
\begin{equation*}
\mathcal{C}_{\mathcal{G}}\leq \zeta(\alpha(x, y)d(Tx, Ty), \beta(M(x, y)) M(x, y)),
\end{equation*}
where
\begin{eqnarray*}
 M(x, y)&=&\max\{d(x, y), d(x, Tx), d(y, Ty)\}\\
 &=& d(x, y).
 \end{eqnarray*}
 This implies
\begin{eqnarray*}
\mathcal{C}_{\mathcal{G}}&\leq& \zeta(\alpha(x, y)d(x, y), \beta(d(x, y)) d(x, y))\\
&<& \mathcal{G}(\beta(d(x, y)) d(x, y)), \alpha(x, y)d(Tx, Ty)),
\end{eqnarray*}
by definition of $\mathcal{G}$, we have
\begin{eqnarray*}
\alpha(x, y)d(x, y)< \beta(d(x, y)) d(x, y).
\end{eqnarray*}
From conditions (2) and (3) and the fact that $\beta \in \mathcal{F}$, we have $\alpha(x, y)\geq1$ and so
\begin{eqnarray*}
d(x, y)\leq \alpha(x, y)d(x, y)< \beta(d(x, y)) d(x, y) <d(x, y),
\end{eqnarray*}
which is a contradiction. Hence $T$ has at most one fixed point.
\end{proof}
\end{lemma}
\begin{lemma}\label{Th4.1}
Let $(X, d)$ be a complete metric space, $\alpha:X\times X\rightarrow [0,
\infty)$ and $T: X \rightarrow X$ be two  functions. Let $\{x_{n}\}$ be a
Picard sequence with initial point at $x_{1}\in X$. Suppose that the
following conditions are satisfied:
\begin{description}
  \item [(1)]  $T$ is Suzuki type $\mathcal{Z}_{(\alpha,
      \mathcal{G})}$-Geraghty contraction;
  \item [(2)] $T$ is triangular $\alpha$-admissible;
  \item [(3)] there exists $x_{1} \in X $ such that $ \alpha (x_{1},
      Tx_{1})\geq 1$.
\end{description}
Then $$\lim\limits_{n\rightarrow\infty}d(x_{n}, x_{n+1})= 0.$$
\begin{proof}
Let $x_{1}\in X$ be such that $\alpha (x_{1}, Tx_{1})\geq 1$ and let
$\{x_{n}\}$ be a Picard sequence, that is, $x_{n}=Tx_{n-1}$, for all $n\in
\mathbb{N}$. If there exist $n\in \mathbb{N}$ such that $d(x_{n}, Tx_{n})= 0$, then $x_{n}$ becomes a fixed
point of $T$, and the proof is complete. Therefore, suppose that
$$0< d(x_{n}, Tx_{n}) ~ ~ \text{for all} ~ ~  n\in\mathbb{N}.$$
Hence, we have
\begin{equation*}
\frac{1}{2}d(x_{n}, Tx_{n})< d(x_{n}, Tx_{n}) = d(x_{n}, x_{n+1}).
\end{equation*}
Since $T$ is  Suzuki type $\mathcal{Z}_{(\alpha, \mathcal{G})}$ -Geraghty
contraction, we have
\begin{eqnarray*}
\mathcal{C}_{\mathcal{G}}&\leq& \zeta(\alpha(x_{n}, x_{n+1})d(Tx_{n}, Tx_{n+1}), \beta(M(x_{n}, x_{n+1})) M(x_{n}, x_{n+1}))\\
&<&\mathcal{G}(\beta(M(x_{n}, x_{n+1})) M(x_{n}, x_{n+1}), \alpha(x_{n}, x_{n+1})d(Tx_{n}, Tx_{n+1})),
\end{eqnarray*}
by definition of $\mathcal{G}$, we get that
\begin{eqnarray*}
 \alpha(x_{n}, x_{n+1})d(x_{n+1}, x_{n+2})< \beta(M(x_{n}, x_{n+1})) M(x_{n}, x_{n+1}),
\end{eqnarray*}
by Lemma \ref{L3.1}, we get
\begin{eqnarray}\label{Eq4.1'}
d(x_{n+1}, x_{n+2})&\leq & \alpha(x_{n}, x_{n+1})d(x_{n+1}, x_{n+2}) \nonumber\\
&<& \beta(M(x_{n}, x_{n+1}))M(x_{n}, x_{n+1}),
\end{eqnarray}
since $\beta \in \mathcal{F}$, we have
\begin{equation}\label{Eq4.1}
d(x_{n+1}, x_{n+2})<M(x_{n}, x_{n+1}),
\end{equation}
where
\begin{eqnarray*}
M(x_{n}, x_{n+1})&\leq& \max\{d(x_{n}, x_{n+1}), d(x_{n}, T x_{n}), d(x_{n+1}, Tx_{n+1})\}\nonumber \\
&=&\max\{d(x_{n}, x_{n+1}), d(x_{n}, x_{n+1}), d(x_{n+1}, x_{n+2})\}\nonumber \\
&=&\max\{d(x_{n}, x_{n+1}), d(x_{n+1}, x_{n+2})\}.
\end{eqnarray*}
If $M(x_{n}, x_{n+1})=d(x_{n+1}, x_{n+2})$ then \eqref{Eq4.1} become
\begin{eqnarray*}
d(x_{n+1}, x_{n+2})&<& d(x_{n+1}, x_{n+2}),
\end{eqnarray*}
which is a contradiction. Hence
\begin{equation}\label{2.4}
M(x_{n}, x_{n+1})=d(x_{n}, x_{n+1}),
\end{equation}
 which
implies that ${d(x_{n}, x_{n+1})}$ is monotonically decreasing
sequence of non negative reals. So there is some $r\geq0$ such that
$$\lim\limits_{n\rightarrow\infty} d(x_{n}, x_{n+1})= r.$$
 We claim that $r=0$, suppose on contrary that $r > 0$, then from (\ref{Eq4.1'}), we have
\begin{eqnarray*}
\frac {d(x_{n+1}, x_{n+2})}{d(x_{n}, x_{n+1})} \leq  \beta (d(x_{n}, x_{n+1}) < 1
\end{eqnarray*}
 which yields that $\lim\limits_{ n\rightarrow \infty }\beta (d(x_{n}, x_{n+1})= 1$.
 since $\beta \in \mathcal{F} $, we derive that
 \begin{equation*}
 \lim\limits_{ n\rightarrow \infty}  d(x_{n}, x_{n+1})= 0,
 \end{equation*}
 which completes the proof.
\end{proof}
\end{lemma}
\begin{lemma}\label{Th4.2}
Let $(X, d)$ be a complete metric space, $\alpha:X\times X\rightarrow [0,
\infty)$ and $T: X \rightarrow X$ be two functions. Suppose that the
following conditions are satisfied:
\begin{description}
  \item [(1)]  $T$ is Suzuki type $\mathcal{Z}_{(\alpha,
      \mathcal{G})}$-Geraghty contraction;
  \item [(2)] $T$ is triangular $\alpha$-admissible;
  \item [(3)] there exists $x_{1} \in X $ such that $ \alpha (x_{1},
      Tx_{1})\geq 1$.
\end{description}
 Then the sequence $\{x_{n}\}$ is bounded.
\begin{proof}
Let $x_{1}\in X$ be such that $\alpha (x_{1}, Tx_{1})\geq 1$. Suppose that
$\{x_{n}\}$ is not bounded sequence. Then there is a subsequence
$\{x_{n_{k}}\}$ of $\{x_{n}\}$ such that $n_{1}=1$ and for each $k\in
\mathbb{N}$, $n_{k+1}$ is the minimum integer such that
 $$ d(x_{n_{k+1}}, x_{n_{k}}) > 1$$ and
$$d(x_{m}, x_{n_{k}}) \leq 1 ~ ~ for ~ ~ n_{k}\leq m\leq n_{k+1}-1.$$
Thus, by triangle inequality, we get
\begin{eqnarray*}
1< d(x_{n_{k+1}}, x_{n_{k}}) &\leq& d(x_{n_{k+1}}, x_{n_{k+1}-1}) +d(x_{n_{k+1}-1}, x_{n_{k}})\\
&\leq& d(x_{n_{k+1}}, x_{n_{k+1}-1})+1.
\end{eqnarray*}
Letting $k \rightarrow \infty$ and using Lemma \ref{Th4.1}, we obtain
$$ \lim\limits_{k\rightarrow \infty} d(x_{n_{k+1}}, x_{n_{k}})= 1.$$
Since
\begin{equation*}
\frac{1}{2}d(x_{n_{k}-1}, Tx_{n_{k}-1})= \frac{1}{2}d(x_{n_{k}-1}, x_{n_{k}})<d(x_{n_{k-1}}, x_{n_{k+1}-1}),
\end{equation*}
then by definition of the mapping $T$, we get that
\begin{eqnarray*}
\mathcal{C}_{\mathcal{G}}&\leq & \zeta(\alpha(d(x_{n_{k+1}-1}, x_{n_{k}-1}))d(x_{n_{k+1}}, x_{n_{k}}), \beta( M(x_{n_{k+1}-1}, x_{n_{k}-1})) M(x_{n_{k+1}-1},
x_{n_{k}-1})) \\
&<& \mathcal{G}(\beta( M(x_{n_{k+1}-1}, x_{n_{k}-1})) M(x_{n_{k+1}-1},
x_{n_{k}-1})), \alpha(d(x_{n_{k+1}-1}, x_{n_{k}-1}))d(x_{n_{k+1}}, x_{n_{k}}))
\end{eqnarray*}
and by definition of $\mathcal{G}$, we have
$$\alpha(d(x_{n_{k+1}-1}, x_{n_{k}-1}))d(x_{n_{k+1}}, x_{n_{k}})< \beta( M(x_{n_{k+1}-1}, x_{n_{k}-1})) M(x_{n_{k+1}-1},
x_{n_{k}-1}).$$
By Lemma \ref{L3.1} and the fact that $\beta\in \mathcal{F}$, we obtain
\begin{eqnarray}\label{Eq4.2}
d(x_{n_{k+1}}, x_{n_{k}})&\leq&\alpha(d(x_{n_{k+1}-1}, x_{n_{k}-1}))d(x_{n_{k+1}}, x_{n_{k}})\nonumber\\
&<& \beta( M(x_{n_{k+1}-1}, x_{n_{k}-1})) M(x_{n_{k+1}-1}, x_{n_{k}-1})\nonumber\\
&<& M(x_{n_{k+1}-1}, x_{n_{k}-1}).
\end{eqnarray}
 Now
\begin{eqnarray}\label{Eq4.3}
1&<& d(x_{n_{k+1}}, x_{n_{k}})<  M(x_{n_{k+1}-1}, x_{n_{k}-1})\nonumber\\
&=&\max\{d(x_{n_{k+1}-1}, x_{n_{k}-1}), d(x_{n_{k+1}-1}, x_{n_{k+1}}),d(x_{n_{k}-1}, x_{n_{k}})\}\nonumber\\
&=&\max\{d(x_{n_{k+1}-1}, x_{n_{k}})+d(x_{n_{k}}, x_{n_{k}-1}), d(x_{n_{k+1}-1}, x_{n_{k+1}}),d(x_{n_{k}-1}, x_{n_{k}})\}\nonumber\\
&\leq&\max\{1+d(x_{n_{k}}, x_{n_{k}-1}), d(x_{n_{k+1}-1}, x_{n_{k+1}}),d(x_{n_{k}-1}, x_{n_{k}})\}.
\end{eqnarray}
Letting $k \rightarrow \infty$, we obtain
\begin{equation*}
1\leq \lim\limits_{k\rightarrow\infty} M(x_{n_{k+1}-1}, x_{n_{k}-1})\leq 1,
\end{equation*}
that is,
\begin{equation*}
\lim\limits _{k\rightarrow\infty}M(x_{n_{k+1}-1}, x_{n_{k}-1})=1.
\end{equation*}
From (\ref{Eq4.2}), we obtain that
\begin{equation*}
\lim\limits _{k\rightarrow\infty}\alpha(d(x_{n_{k+1}-1}, x_{n_{k}-1}))d(x_{n_{k+1}}, x_{n_{k}})=1,
\end{equation*}
and
\begin{equation*}
\lim\limits _{k\rightarrow\infty}\beta( M(x_{n_{k+1}-1}, x_{n_{k}-1})) M(x_{n_{k+1}-1}, x_{n_{k}-1})=1.
\end{equation*}
Further, since $$\frac{1}{2}d(x_{n_{k}-1}, Tx_{n_{k}-1}) < d(x_{n_{k}-1},
x_{n_{k+1}-1}),$$  we have by definition of $T$
\begin{eqnarray*}
\mathcal{C}_{\mathcal{G}}&\leq& \zeta(\alpha(x_{n_{k+1}-1}, x_{n_{k}-1})d(Tx_{n_{k+1}-1}, Tx_{n_{k}-1}), \beta(M(x_{n_{k+1}-1}, x_{n_{k}-1}))M(x_{n_{k+1}-1}, x_{n_{k}-1}))\\
&\leq& \lim\limits_{k\rightarrow\infty}\sup\zeta(\alpha(x_{n_{k+1}-1}, x_{n_{k}-1})d(x_{n_{k+1}}, x_{n_{k}}), \beta(M(x_{n_{k+1}-1}, x_{n_{k}-1}))M(x_{n_{k+1}-1}, x_{n_{k}-1})) \\ &<&\mathcal{C}_{\mathcal{G}}
\end{eqnarray*}
a contradiction. Hence $\{x_{n}\}$ is a bounded sequence.
\end{proof}
\end{lemma}
\begin{lemma}\label{L4.3}
 $(X, d)$ be a complete metric space, $\alpha:X\times X\rightarrow [0,
\infty)$ and $T: X \rightarrow X$ be two functions. Suppose that the
following conditions are satisfied:
\begin{description}
  \item [(1)]  $T$ is  Suzuki type $\mathcal{Z}_{(\alpha,
      \mathcal{G})}$-Geraghty contraction;
  \item [(2)] $T$ is triangular $\alpha$-admissible;
  \item [(3)] there exists $x_{1} \in X $ such that $ \alpha (x_{1},
      Tx_{1})\geq 1$.
\end{description}
Then the sequence $\{x_{n}\}$ is Cauchy.
\begin{proof}
To show that $\{x_{n}\}$ is a Cauchy sequence, we suppose on contrary that $\{x_{n}\}$ it is not. Then by Lemma \ref{L4.2}, we have
$$\lim\limits_{k\rightarrow\infty}d(x_{m(k)},x_{n(k)})=\lim\limits_{k\rightarrow\infty} d(x_{m(k)+1},x_{n(k)+1})=\varepsilon$$
and consequently
$$\lim\limits_{k\rightarrow\infty} M(x_{m(k)},x_{n(k)})=\varepsilon.$$
By Lemma \ref{L4.2} and \ref{Th4.2}, we have
\begin{eqnarray*}
\frac{1}{2}d(x_{m(k)}, Tx_{m(k)})=\frac{1}{2}d(x_{m(k)}, x_{m(k)+1})<d(x_{m(k)}, x_{n(k)}),
\end{eqnarray*}
 then by definition of $T$, we have
\begin{eqnarray*}
\mathcal{C}_{\mathcal{G}} &\leq& \zeta(\alpha(x_{m(k)}, x_{n(k)})d(Tx_{m(k)}, Tx_{n(k)}), \beta(M(x_{m(k)}, x_{n(k)}))M(x_{m(k)}, x_{n(k)}))\\
&<&\mathcal{G}(\beta(M(x_{m(k)}, x_{n(k)}))M(x_{m(k)}, x_{n(k)}), \alpha(x_{m(k)}, x_{n(k)})d(Tx_{m(k)}, Tx_{n(k)})).
\end{eqnarray*}
By definition of ${\mathcal{G}}$, we have
\begin{eqnarray*}
\alpha(x_{m(k)}, x_{n(k)})d(Tx_{m(k)}, Tx_{n(k)})<\beta(M(x_{m(k)}, x_{n(k)}))M(x_{m(k)}, x_{n(k)}).
\end{eqnarray*}
Since $T$ is triangular $\alpha$-admissible, so by Lemma \ref{L3.1}, we have
$\alpha(x_{m(k)},x_{n(k)})\geq1$ and $\beta\in \mathcal{F}$, we have
\begin{eqnarray}\label{Eq4.7'}
d(Tx_{m(k)}, Tx_{n(k)})&<&\alpha(x_{m(k)}, x_{n(k)})d(Tx_{m(k)}, Tx_{n(k)}) \nonumber \\
&<&\beta(M(x_{m(k)}, x_{n(k)}))M(x_{m(k)}, x_{n(k)}) \nonumber\\
&<& M(x_{m(k)}, x_{n(k)}),
\end{eqnarray}
by \eqref{Eq4.7'}, we obtain
\begin{equation*}
\lim\limits _{k\rightarrow\infty}\alpha(x_{m(k)}, x_{n(k)})d(Tx_{m(k)}, Tx_{n(k)})=\varepsilon,
\end{equation*}
and
\begin{equation*}
\lim\limits _{k\rightarrow\infty}\beta(M(x_{m(k)}, x_{n(k)}))M(x_{m(k)}, x_{n(k)})=\varepsilon.
\end{equation*}
So by definition of $\mathcal{C}_{\mathcal{G}}$ simulation function, we have
\begin{eqnarray*}
\mathcal{C}_{\mathcal{G}} &\leq& \zeta(\alpha(x_{m(k)}, x_{n(k)})d(Tx_{m(k)}, Tx_{n(k)}), \beta(M(x_{m(k)}, x_{n(k)}))M(x_{m(k)}, x_{n(k)}))\\
&\leq&\lim\limits_{k\rightarrow\infty}\sup\zeta(\alpha(x_{m(k)}, x_{n(k)})d(Tx_{m(k)}, Tx_{n(k)}), \beta(M(x_{m(k)}, x_{n(k)}))M(x_{m(k)}, x_{n(k)}))\\
&<&\mathcal{C}_{\mathcal{G}}
\end{eqnarray*}
a contradiction. Hence $\{x_{n}\}$ is a Cauchy sequence.
\end{proof}
\end{lemma}
We now prove the existence result.
\begin{theorem}\label{Th4.4}
 $(X, d)$ be a complete metric space, $\alpha:X\times X\rightarrow [0,
\infty)$ and $T: X \rightarrow X$ be two functions. Suppose that the
following conditions are satisfied:
\begin{description}
  \item [(1)]  $T$ is  Suzuki type $\mathcal{Z}_{(\alpha,
      \mathcal{G})}$-Geraghty contraction ;
  \item [(2)] $T$ is triangular $\alpha$-admissible;
  \item [(3)] there exists $x_{1} \in X $ such that $ \alpha (x_{1},
      Tx_{1})\geq 1$.
\end{description}
Then $T$ has a fixed point.
\begin{proof}
From Lemma \ref{L4.3}, $\{x_{n}\}$ is a Cauchy sequence and by the completeness of $X$ there exists $x\in X$ such that
\begin{eqnarray}\label{Eq4.7}
\lim\limits_{n\rightarrow \infty}x_{n}=x.
\end{eqnarray}
Now, we show that $x$ is a fixed point of $T$. We claim that
\begin{eqnarray*}
\frac{1}{2}d(x_{n}, Tx_{n})<d(x_{n}, x)~ ~ or ~ ~ \frac{1}{2}d(x_{n+1}, Tx_{n+1})<d(x_{n+1}, x), ~ ~for ~ ~ all ~ ~ n\in\mathbb{N}.
\end{eqnarray*}
That is,
\begin{eqnarray}\label{Eq4.8}
(I):=~\frac{1}{2}d(x_{n}, Tx_{n})<d(x_{n}, x) ~ ~ \text{or} ~ ~ (II):=~\frac{1}{2}d(Tx_{n},T^{2}x_{n})<d(Tx_{n}, x), ~ ~ for ~ ~ all~ ~ n\in\mathbb{N}.
\end{eqnarray}
Assume that there exists $m\in\mathbb{N}$ such that
\begin{eqnarray}\label{Eq4.9}
\frac{1}{2}d(x_{m}, Tx_{m})\geq d(x_{m}, x) ~ ~ or ~ ~ \frac{1}{2}d(Tx_{m},T^{2}x_{m})\geq d(Tx_{m}, x).
\end{eqnarray}
Hence,
\begin{equation*}
2d(x_{m}, x)\leq d(x_{m},Tx_{m})\leq d(x_{m}, x)+d(x, Tx_{m})
\end{equation*}
this implies that
\begin{eqnarray}\label{Eq4.10}
d(x_{m}, x)\leq d(x, Tx_{m}).
\end{eqnarray}
From (\ref{2.4})and (\ref{Eq4.10}), we have
\begin{eqnarray}\label{Eq4.11}
d(Tx_{m},T^{2}x_{m})< d(x_{m},Tx_{m})\leq d(x_{m}, x)+d(x, Tx_{m})< 2d(x, Tx_{m}).
\end{eqnarray}
It follows from (\ref{Eq4.9}) and (\ref{Eq4.11}), that $d(Tx_{m},T^{2}x_{m})<
d(Tx_{m},T^{2}x_{m})$, which is a contradiction. Thus (\ref{Eq4.8}) holds. If
part (I) of (\ref{Eq4.8}) is true, by definition of $T$, we have
\begin{eqnarray*}
\mathcal{C}_{\mathcal{G}} &\leq& \zeta(\alpha(x_{n}, x)d(Tx_{n}, Tx), \beta(M(x_{n}, x))M(x_{n}, x))\\
&<& \mathcal{G}(\beta(M(x_{n}, x))M(x_{n}, x), \alpha(x_{n}, x)d(Tx_{n}, Tx)).
\end{eqnarray*}
By definition of $\mathcal{G}$, we have
\begin{eqnarray*}
\alpha(x_{n}, x)d(Tx_{n}, Tx)<\beta(M(x_{n}, x))M(x_{n}, x),
\end{eqnarray*}
by Lemma \ref{L3.1} and the fact that $\beta \in \mathcal{F}$, we have
\begin{eqnarray*}
d(Tx_{n}, Tx)&\leq & \alpha(x_{n}, x)d(Tx_{n}, Tx)\\
&<& \beta(M(x_{n}, x))M(x_{n}, x)\\
&<& M(x_{n}, x).
\end{eqnarray*}
This implies
\begin{eqnarray*}
d(Tx_{n}, Tx)&<& M(x_{n}, x)= \max\{d(x_{n}, x), d(x_{n}, Tx_{n}), d(x, Tx)\}.
\end{eqnarray*}
Letting $n\rightarrow \infty$ and using (\ref{Eq4.7}), we obtain
\begin{eqnarray}\label{Eq4.12}
\lim\limits_{n\rightarrow \infty}M(x_{n}, x)=d(x, Tx).
\end{eqnarray}
Now,
\begin{eqnarray*}
\mathcal{C}_{\mathcal{G}}&\leq& \zeta(\alpha(x_{n}, x)d(Tx_{n}, Tx), \beta(M(x_{n}, x))M(x_{n}, x))\\
&\leq& \lim\limits_{n\rightarrow\infty}\sup\beta(M(x_{n}, x))d(x, Tx)-\alpha(x_{n}, x)d(x, Tx)\\
&<&\mathcal{C}_{\mathcal{G}}
\end{eqnarray*}
a contradiction. Hence $x=Tx$, {it i.e.}  $x$ is a fixed point of $T$.
 On the other hand if part (II) of (\ref{Eq4.8}) is true, again by definition of $T$, we have
\begin{eqnarray*}
\mathcal{C}_{\mathcal{G}}&\leq& \zeta(\alpha(Tx_{n}, x)d(T^{2}x_{n}, Tx), \beta(M(Tx_{n}, x))M(Tx_{n}, x))\\
&<& \mathcal{G}(\beta(M(Tx_{n}, x))M(Tx_{n}, x), \alpha(Tx_{n}, x)d(T^{2}x_{n}, Tx)),
\end{eqnarray*}
by Lemma \ref{L3.1} and the fact that $\beta \in \mathcal{F}$, we have
\begin{eqnarray*}
d(T^{2}x_{n}, Tx)&\leq &\alpha(Tx_{n}, x)d(T^{2}x_{n}, Tx)\\
&<& \beta(M(Tx_{n}, x))M(Tx_{n}, x)\\
&<& M(Tx_{n}, x)
\end{eqnarray*}
this implies
\begin{eqnarray*}
d(T^{2}x_{n}, Tx)&<& M(Tx_{n}, x)= \max\{d(Tx_{n}, x), d(Tx_{n}, T^{2}x_{n}), d(x, Tx)\}.
\end{eqnarray*}
Letting $n\rightarrow \infty$ and by using (\ref{Eq4.7}), we obtain
\begin{eqnarray}\label{Eq4.13}
\lim\limits_{n\rightarrow \infty}M(Tx_{n}, x)=d(x, Tx).
\end{eqnarray}
Now,
\begin{eqnarray*}
\mathcal{C}_{\mathcal{G}}&\leq& \zeta(\alpha(Tx_{n}, x)d(T^{2}x_{n}, Tx), \beta(M(Tx_{n}, x))M(Tx_{n}, x))\\
&\leq& \lim\limits_{n\rightarrow\infty}\sup\beta(M(Tx_{n}, x))d(x, Tx)-\alpha(Tx_{n}, x)d(x, Tx)\\
&<&\mathcal{C}_{\mathcal{G}}
\end{eqnarray*}
a contradiction. Hence $x=Tx$, {\it i.e.} $x$ is a fixed point of $T$.
Uniqueness of fixed point follows from Lemma \ref{L4.1}.
\end{proof}
\end{theorem}
\begin{example}\label{Ex4.1}
Let $X=\{0, \frac{1}{2}, 1, 2\}$ and $d:X \times X\rightarrow R$ be defined
by $d(x, y)= |x- y|$ for all $x, y \in X$. Let $\zeta(t, s)= \frac{s}{s+1} -
t$, $\mathcal{G}(s, t)=s-t$ for all $t, s\in[0, \infty)$ , $\mathcal{C}(\mathcal{G})=0$
and $\beta(t)=\frac{1}{1+t}$ for all $t\geq 0$. It is clear that $\beta
\in \mathcal{F}$. Define $T: X\rightarrow X$ by
$$Tx= \begin{cases}
\frac {1}{2}    & ~ if ~ x\in\{0,\frac{1}{2}, 1\}, \\
 \frac {1}{9}    & ~ if ~ x\in \{2\}.
\end{cases}$$
and a function $\alpha: X\times X \rightarrow [0, \infty)$ by
$$\alpha(x, y) = \begin{cases}
1 & ~ if ~ x, y\in \{0, \frac{1}{2}, 1\} \\
0 & otherwise.
\end{cases}$$
Condition $(3)$ of Theorem \ref{Th4.4} is satisfied with $x_{1}= 1$. Let $x, y
\in X$ be such that $\alpha(x, y)\geq 1$. Then $x, y \in \{0,\frac{1}{2}, 1\}$,
and so $Tx=\frac{1}{2}=Ty$ and $\alpha(Tx, Ty)= 1$. Hence $T$ is
$\alpha$- admissible, hence $(2)$ is satisfied. We show that condition $(1)$
of theorem $(4.4)$ is satisfied. For $ x, y\in X$, with $x<y$ we have
\begin{eqnarray*}
\frac{1}{2}d(x, Tx)<d(x, y),
\end{eqnarray*}
and
\begin{eqnarray*}
\zeta(\alpha(x, y)d(Tx, Ty), \beta(M(x, y))M(x, y))= \frac{\beta(M(x, y))M(x, y)}{\beta(M(x, y))M(x, y)+1}-\alpha(x, y)d(Tx, Ty),
\end{eqnarray*}
where
\begin{eqnarray*}
M(x, y)&=&\max\{d(x, y), d(x, Tx), d(y, Ty)\}.
\end{eqnarray*}
Now,
\begin{eqnarray*}
\zeta(\alpha(x, y)d(Tx, Ty), \beta(M(x, y))M(x, y))&=& \frac{\beta(M(x, y))M(x, y)}{\beta(M(x, y))M(x, y)+1}-\alpha(x, y)d(Tx, Ty)\\
&=& \frac{\frac{M(x, y)}{M(x, y)+1}}{\frac{M(x, y)}{M(x, y)+1}+1}-d(Tx, Ty)\\
&=&\frac{M(x, y)}{1+2M(x, y)}>0.
\end{eqnarray*}
Hence, $T$ is  Suzuki type $\mathcal{Z}_{(\alpha, \mathcal{G})}$-Geraghty contraction with respect to
$\zeta$. Thus all the conditions of Theorem \ref{Th4.4} are satisfied and $T$
has a unique fixed point $x=\frac{1}{2}\in X$. Since $T$ is not continuous,
then it is not $\mathcal{Z}$-contraction and so not contractive.
\end{example}
\section{Application}
In this section, we present an application of Theorem \ref{Th4.4} to guarantee the
existence and uniqueness of solution of the following non linear
Hammerstein integral equation:
\begin{equation}\label{Eq5.1}
x(t)= f(t)+\int_{0}^{t} K(t, s)h(s, x(s))ds,
\end{equation}
where the unknown function $x(t)$ takes real values. Let $X= C([0, 1])$ be
the space of all real continuous functions defined on $[0, 1]$. It is well
known that $C[0, 1]$ endowed with the metric
\begin{equation}\label{Eq5.2}
d(x, y)= \|x-y\|= \max_{t\in[0, 1]}|x(t)-y(t)|.
\end{equation}
So, $(X, d)$ is a complete metric space. Define a mapping $T: X\rightarrow X$
by
\begin{equation}\label{Eq5.3}
T(x)(t)= f(t)+\int_{0}^{t} K(t, s)h(s, x(s))ds,
\end{equation}
for all $t\in[0, 1]$.

Suppose that the following conditions are satisfied:
\begin{description}
  \item [(1)] $f\in C([0, 1])\times (-\infty, +\infty), f\in X$ and $K\in
      C([0, 1]\times[0, 1])$ such that $K(t, s)\geq0$;
  \item [(2)] $h(t, s):(-\infty, +\infty))\rightarrow(-\infty, +\infty))
      $ is increasing for all $t\in[0, 1]$, such that
      $$ \frac{1}{2}d(x, F(x))< d(x, y)$$
      implies $|h(t, x)- h(t, y)|< M(x, y)$, for all distinct $x, y\in
  X,t\in[0, 1] $, where
  $$M(x, y)= \max\{|x-y|, |x-Tx|, |y-Ty|\};$$
  \item [(3)] for all $t,s\in [0, 1]$, $$\int_{0}^{t}|K(t, s)|ds \leq
      \frac{1}{1+M(x, y)}.$$
\end{description}
\begin{theorem}\label{Th5.1}
Under the assumptions (1)-(3), the integral equation \eqref{Eq5.1} has a unique solution in $X=C([0, 1])$.

\begin{proof}
Define a mapping $\alpha:C[0,1]\times C[0,1]\to \mathbb{R}$ by
 \begin{eqnarray*}
\alpha(x,y)=\left\{\begin{array}{cc}
                               1 & \hbox{if}~~x(s),y(s)\in[0,\infty)~\text{for all}~~ s\in[0,1] \\
                               0 & \text{otherwise}.
                             \end{array}
                             \right.
\end{eqnarray*}
Notice that the existence of solution to \eqref{Eq5.1} is equivalent to the existence of fixed point of $T$. Now, we will show that all hypothesis of Theorem \ref{Th4.4} are satisfied.
We first show that $T$ is triangular $\alpha$-admissible mapping. Indeed, for $x,y, z\in C[0, 1]$ such that $\alpha(x(s),y(s))\geq1$ and $\alpha(y(s),z(s))\geq1$, we have $x(s),y(s), z(s)\geq0$ for all $s\in [0,1]$. Therefore $\alpha(x(s),z(s))\geq1$ and hence $T$ is triangular $\alpha$-admissible mapping. Also there is $x_{1}(s)\geq 0$ such that $Tx_{1}(s) \geq 0$, then $\alpha(x_{1}(s), Tx_{1}(s))\geq 1$. \\
Now, we claim that the mapping $T:X\rightarrow X$ define by \eqref{Eq5.3} is
a Suzuki type $\mathcal{Z}_{(\alpha, \mathcal{G})}$ Geraghty contraction.

From condition (2) and (3), for all distinct $x, y\in C([0, 1]), t\in[0, 1]$,
we have
\begin{eqnarray*}
\alpha(x, y)|Tx(t)-Ty(t)|&=&|Tx(t)-Ty(t)| \\
 &=& \Bigg|\int_{0}^{t} K(t, s)(h(s, x(s))-h(s, y(s)))ds \Bigg|\\
&\leq& \int_{0}^{t} |K(t, s)||h(s, x(s))-h(s, y(s))|ds\\
&<&\int_{0}^{t}|K(t, s)|M(x(s), y(s))ds,
\end{eqnarray*}
where
\begin{eqnarray*}
M(x(s), y(s))&=&\max\{|x(s)-y(s)|, |x(s)-Tx(s)|, |y(s)-Ty(s)|\}\\
&\leq& \max\{d(x, y), d(x, Tx), d(y, Ty)\}\\
&=&M(x, y).
\end{eqnarray*}
This implies,
\begin{eqnarray*}
\alpha(x, y)|Tx(t)-Ty(t)|&<&M(x, y)\int_{0}^{t}|K(t, s)|ds \\
&\leq &\frac{M(x, y)}{1+M(x, y)}.
\end{eqnarray*}
Hence, the mapping $T$ is a Suzuki type $\mathcal{Z}_{(\alpha, \mathcal{G})}$
 Geraghty contraction with $\zeta(t,s)=\mathcal{G}(s, t)=s-t$ and $\beta(t)=\frac{1}{1+t}$. Thus we can apply Theorem \ref{Th4.4}, which guarantee
 the existence of a unique fixed point $w\in X$. That is, $w$ is the unique
 solution of non linear Hammerstein integral equation \eqref{Eq5.1}.
\end{proof}
\end{theorem}

\end{document}